# Symplectification of Circular Arcs and Arc Splines

*An Efficient Representation of Planar Circular Arcs for Use in Polyarcs and Arc Splines*

STEFAN GÖSSNER[1]

[1]Dortmund University of Applied Sciences. Department of Mechanical Engineering

## *Abstract*

In this article, circular arcs are considered both individually and as elements of a piecewise circular curve. The endpoint parameterization proves to be quite advantageous here. The perspective of symplectic geometry provides new vectorial relationships for the circular arc. Curves are considered whose neighboring circular elements each have a common end point or, in addition, a common tangent. These arc splines prove to be a one-parameter curve family, whereby this parameter can be optimized with regard to various criteria.

# Content

# 1. Introduction

Circular arcs play a fundamental role in Engineering, Manufacturing and Computer Graphics. Circular arcs are used to define the shapes of mechanical parts during design. They are important in CNC Machining and Robot path planning, where circular and linear interpolation are common [4-9]. Equidistant curves to piecewise circular curves are determined easily.

In this paper circular arcs are considered purely geometrically from the perspective of symplectic geometry [1-3]. The paper is structured as follows. First, a single circular arc and its various parameter sets are discussed. Then, polyarcs and arc splines are presented as planar curves composed of individual arcs. Finally, arc splines are identified as a one parameter family of curves that belong to their underlying polygonial curve. An optimal parameter can be found by minimizing certain scalar curve properties.

## *1.1 Prerequisites*

The *symplectic vector space* extends the Euclidean vector space in $\mathbb{R}^2$ by a *complex structure* $\mathbf{J} = \begin{pmatrix} 0 & -1 \\ 1 & 0 \end{pmatrix}$. This skew-symmetric matrix is also referred to as the *standard symplectic structure*. It is an *orthogonal operator* in $\mathbb{R}^2$ which, with an anti-clockwise rotation by $\pi/2$, transforms any vector into a *skew-orthogonal* complementary one [1,2].

$$\tilde{\mathbf{a}} = \mathbf{J} \begin{pmatrix} a_1 \\ a_2 \end{pmatrix} = \begin{pmatrix} -a_2 \\ a_1 \end{pmatrix}$$

As a shortcut we will place a *tilde* '~' symbol over the skew-orthogonal vector variable, visually reflecting the skew-symmetry of its operator matrix. Taking the *standard scalar product* on the plane [3]

$$\mathbf{a}\mathbf{b} = a_1 b_1 + a_2 b_2$$

and applying the orthogonal operator to the first vector gets us to the *skew-scalar product*

$$\tilde{\mathbf{a}}\mathbf{b} = a_1 b_2 - a_2 b_1 \,.$$

The skew-scalar product gives us the *area* of the parallelogram spanned by two vectors $\mathbf{a}$ and $\mathbf{b}$. We get a *directed* or *oriented* area due to the inherent antisymmetry $\tilde{\mathbf{a}}\mathbf{b} = -\mathbf{a}\tilde{\mathbf{b}}$ [2,3].

Rotating a vector $\mathbf{a}$ by $\alpha$ into vector $\mathbf{a}_\alpha$ gives [3]

$$\mathbf{a}_\alpha = \begin{pmatrix} \cos\alpha & -\sin\alpha \\ \sin\alpha & \cos\alpha \end{pmatrix} \begin{pmatrix} a_1 \\ a_2 \end{pmatrix} = \cos\alpha \begin{pmatrix} a_1 \\ a_2 \end{pmatrix} + \sin\alpha \begin{pmatrix} -a_2 \\ a_1 \end{pmatrix} = \cos\alpha\,\mathbf{a} + \sin\alpha\,\tilde{\mathbf{a}}$$

Many publications dealing with plane vectors agree on a special $\mathbb{R}^2$ form of the cross product [5,7,9].

$$\mathbf{a} \times \mathbf{b} = a_1 b_2 - a_2 b_1 \,.$$

As a result, planar vector equations arise that can hardly be treated further vectorially, but rather have to be broken down into their scalar components. Complex numbers also often suffer from a similar problem.

The symplectic geometry can be advantageously restricted to the standard scalar product. With its fundamental rules [3], vectors in equations remain intact during further treatment. As a result, new equations and formulas are often obtained.

# 2. Single Circular Arc

A circular arc is uniquely defined by a set of five scalar parameters. We would like to discuss two different parameter sets:

- Center based parameters (Fig. 1a)
- Endpoint based parameters (Fig. 1b)

The classic center based parameters are briefly discussed, whereas the endpoint based parameters are illustrated in more detail subsequently.

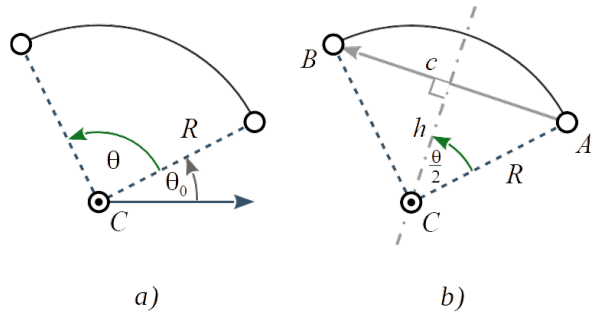

Fig. 1: Different Arc Parameters Sets

## 2.1 Center Based Parameter Set

There is a set of five scalar parameters that uniquely define a circular arc. The classic center based parameters (Fig. 1a) are

- $x_C, y_C$   Center coordinates
- $R$   Radius
- $\theta_0, \theta$   Start angle, angular range $(-2\pi, 2\pi)$

Sometimes start and end angles $\theta_0, \theta_1$ are used instead of the latter two. That parameter set turns out to be ambiguous, as it leads to two complementary arcs running in opposite directions.

## 2.2 Endpoint Parameters

If the circular arc is part of a piecewise circular curve, i.e. polyarc or arc spline, then using endpoint parameters is preferable.

- $x_A, y_A$   Start point coordinates
- $x_B, y_B$   End point coordinates
- $\theta$   Angular range

The circular arc is built over the endpoint (chord) vector $AB$ (Fig. 1b). Four parameters are occupied by both endpoint coordinates. Angle $\theta$ is preferred over the – possibly very big – radius $R$ as the remaining fifth parameter.

Its parameter equation and important properties like length, segment area and bending energy are determined below.

> **Lemma 1**:
> Let $c$ be the distance from an arc's startpoint $A$ to its endpoint $B$. The locus of its center is then necessarily somewhere on the perpendicular bisector of $AB$ and can be controlled either by the angular range $\theta$ or by the radius $R$. The relation between the two is
>
> $$R \sin \frac{\theta}{2} = \frac{c}{2} \quad for \quad c > 0 \qquad (1)$$
>
> A consequence of this is a signed angle $\theta \in (-2\pi, 2\pi)$ and a signed radius, with $|R| \in \left[\frac{1}{2}c, \infty\right)$. Both values always have equal signs constrained by equation (1).

*Proof.* By trigonometric relation in the right-angled half of the isosceles triangle $ABC$ in Fig. 1b. □

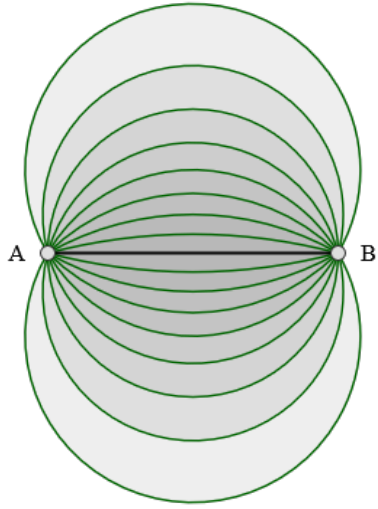

Fig. 2: One parameter family of arcs

The circular arcs in Figure 2 are controlled by angle $\theta$ in the range from -270° to 270° in 30° steps.

The following circular arc formulae solely depend on endpoint parameters, i.e. segment vector $\mathbf{c}$ and angular range $\theta$. In this case, we dispense with knowledge of the absolute position of the arc and content ourselves with three scalar parameters.

## 2.3 Center Vector

> **Lemma 2**:
> Let $\theta$ be an arc's enclosed angle and $\mathbf{c}$ be the vector from its startpoint $A$ to its endpoint $B$. The radiusvector $\mathbf{r}_0$ from the arc's center $C$ to its start point $A$ is then
>
> $$\mathbf{r}_0 = -\frac{\sin \frac{\theta}{2} \mathbf{c} + \cos \frac{\theta}{2} \tilde{\mathbf{c}}}{2 \sin \frac{\theta}{2}} . \qquad (2)$$

*Proof.* The vector loop closure equation of the right-angled half of the isosceles triangle $ABC$ in Fig. 1b reads

$$\frac{\mathbf{c}}{2} + \frac{\tilde{\mathbf{c}}}{2 \tan \frac{\theta}{2}} + \mathbf{r}_0 = \mathbf{0} .$$

Resolving for $\mathbf{r}_0$, using $\tan x = \frac{\sin x}{\cos x}$ and bringing the right side to a common denominator leads to symmetric equation (2). □

## 2.4 Parametric Vector Equation

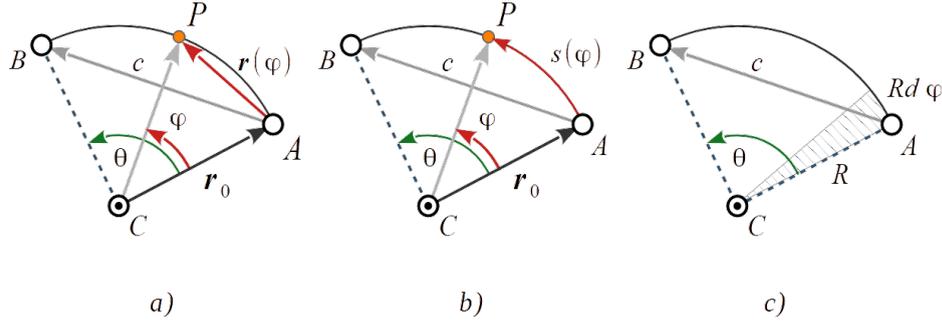

Fig. 2: Parametrizing arc, arc length, arc segment area

> **Theorem 1**:
> Let $\theta$ be an arc's enclosed angle and $\mathbf{c}$ be the vector from its startpoint to its endpoint. The parametric vector equation of the arc with respect to its startpoint is
>
> $$\mathbf{r}(u) = \begin{cases} \dfrac{\sin \frac{u\theta}{2}}{\sin \frac{\theta}{2}} \left[ \cos \frac{(1-u)\theta}{2} \mathbf{c} - \sin \frac{(1-u)\theta}{2} \tilde{\mathbf{c}} \right] & \text{if } \theta \neq 0 \\ u\mathbf{c} & \text{if } \theta = 0 \end{cases} \quad for \quad u \in [0,1] \qquad (3)$$

*Proof.* The loop closure equation of the triangle $ACP$ in Fig. 2a gives the parametric equation of the circular arc

$$\mathbf{r}(\varphi) = -\mathbf{r}_0 + \cos\varphi\, \mathbf{r}_0 + \sin\varphi\, \tilde{\mathbf{r}}_0 = (\cos\varphi - 1)\mathbf{r}_0 + \sin\varphi\, \tilde{\mathbf{r}}_0 \,. \qquad (3.1)$$

Using the half-angle rules $1 - \cos x = 2\sin^2 \frac{x}{2}$ and $\sin x = 2 \sin \frac{x}{2} \cos \frac{x}{2}$ yields

$$\mathbf{r}(\varphi) = 2 \sin \frac{\varphi}{2} \left( -\sin \frac{\varphi}{2} \mathbf{r}_0 + \cos \frac{\varphi}{2} \tilde{\mathbf{r}}_0 \right) . \qquad (3.2)$$

Inserting the radiusvector $\mathbf{r}_0$ from (2) leads to

$$\mathbf{r}(\varphi) = \frac{\sin \frac{\varphi}{2}}{\sin \frac{\theta}{2}} \left[ \left( \sin \tfrac{\varphi}{2} \sin \tfrac{\theta}{2} + \cos \tfrac{\varphi}{2} \cos \tfrac{\theta}{2} \right) \mathbf{c} + \left( \sin \tfrac{\varphi}{2} \cos \tfrac{\theta}{2} - \cos \tfrac{\varphi}{2} \sin \tfrac{\theta}{2} \right) \tilde{\mathbf{c}} \right] ,$$

which can be reduced by addition theorems to

$$\mathbf{r}(\varphi) = \frac{\sin \frac{\varphi}{2}}{\sin \frac{\theta}{2}} \left( \cos \tfrac{\theta - \varphi}{2} \mathbf{c} - \sin \tfrac{\theta - \varphi}{2} \tilde{\mathbf{c}} \right) . \qquad (3.3)$$

Angle $\varphi \in [0, \theta]$ might be expressed favourably by $\varphi = u\theta$, $u \in [0,1]$, which avoids the need to care about synchronizing the signs of $\varphi$ and $\theta$. This leads to the non-zero case in equation (1).

We are allowed to write for small angles $\sin x \approx x$ and $\cos x \approx 1$, which gives us for a small $\theta$

$$\mathbf{r}(u) = u \left( \mathbf{c} - \tfrac{(1-u)\theta}{2} \tilde{\mathbf{c}} \right) .$$

A vanishing $\theta$ herein results in $\mathbf{r}(u) = u\mathbf{c}$ – the zero case of equation (3), i.e. an arc degenerating to a line. □

## 2.5 Arc Length

> **Theorem 2**:
> Let $\theta$ be an arc's enclosed angle and $c$ be the distance from its startpoint to its endpoint. The arc length $s$ parameterized by normalized value $u$ is
>
> $$s(u) = \begin{cases} uc \dfrac{\theta}{2\sin\frac{\theta}{2}} & \text{if } \theta \neq 0 \\ uc & \text{if } \theta = 0 \end{cases} \quad for \quad u \in [0,1] \qquad (4)$$
>
> This equation is arc length parameterized.

*Proof.* The arc length of an arbitrary planar curve is defined by

$$s(\varphi) = \int_0^\varphi \sqrt{\mathbf{r}'(\phi)^2}\, d\phi. \qquad (4.1)$$

Taking the derivation of the circular arc's parametric equation (3.1) with respect to $\varphi$

$$\mathbf{r}'(\varphi) = -\sin\varphi\, \mathbf{r}_0 + \cos\varphi\, \tilde{\mathbf{r}}_0.$$

and inserting gives with

$$s(\varphi) = \int_0^\varphi \sqrt{(\sin^2\phi + \cos^2\phi)R^2}\, d\phi = \varphi R \qquad (4.2)$$

the well known length formula of circular arcs. Again substituting $\varphi = u\theta$ and taking $R$ from (1) leads to

$$s(u) = uc\, \frac{\theta}{2\sin\frac{\theta}{2}}.$$

Angle $\theta \in (-2\pi, 2\pi)$ and value of $\sin\frac{\theta}{2}$ are always of same sign, thus $s(u)$ is always positive. Setting $\sin x \approx x$ for small angles simplifies to linear function $s(u) = uc$ (Fig. 2b).

Equal parameter intervals $\Delta u$ correspond to equal arc lengths $\Delta s$ in equation (4), which thus is arc length parametrized. $\square$

## 2.6 Arc Segment Area

The arc sector area is derived from the small hatched area in Fig. 2c by integration. Subtracting triangle area $\Delta CAB$ from the sector area yields the arc segment area.

> **Theorem 3**:
> Let $\theta$ be an arc's enclosed angle and $c$ be the distance from its startpoint to its endpoint. The signed arc segment area is then
>
> $$A_{seg} = \frac{c^2}{4} \cdot \frac{\theta - \sin\theta}{1 - \cos\theta} \qquad (5)$$

*Proof.* The small hatched arc sector area in Fig. 2c is approximated as a right triangle area $\frac{1}{2}R^2 d\phi$. From this we yield the total sector area

$$A_{sector} = \int_0^\theta \tfrac{1}{2}R^2 d\phi = \tfrac{1}{2}R^2\theta \tag{5.1}$$

The triangle area in Fig. 1b is $A_{CAB} = \tfrac{1}{2}ch$, where $h = R\cos\tfrac{\theta}{2}$. Using expression (1) for $c$ gives

$$A_{CAB} = R^2 \cos\tfrac{\theta}{2} \sin\tfrac{\theta}{2} = \tfrac{1}{2}R^2 \sin\theta\,,$$

while applying the half-angle rule $\sin x = 2\sin\tfrac{x}{2}\cos\tfrac{x}{2}$. Now taking $A_{seg} = A_{sector} - A_{CAB}$ yields

$$A_{seg} = \tfrac{1}{2}R^2(\theta - \sin\theta)\,. \tag{5.2}$$

Equation (5) results with $R^2 = \dfrac{c^2}{4\sin^2\tfrac{\theta}{2}}$ from equation (1) and half-angle term $\sin^2\tfrac{x}{2} = \dfrac{1-\cos x}{2}$. □

The segment area disappears for angle $\theta = 0$. This can be proved by applying L'Hospital's rule twice to expression (5).

## 2.6 Bending Energy

Imagine the circular arc as a formerly straight, elastic beam of homogeneous material and cross-section, that is loaded by a pure bending moment (Fig. 3).

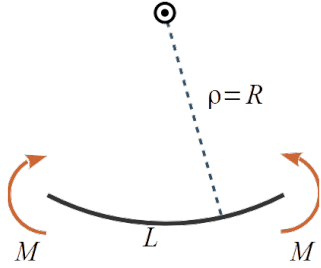

Fig. 3: Elastic beam loaded by a pure moment.

The normal stress distribution in the beam's cross-section along the (in Fig. 3 vertical) y-direction is $\sigma(y) = \dfrac{M}{I}y$ due to load and $\sigma(y) = \dfrac{E}{\rho}y$ due to strain, following the Bernoulli-Euler beam theory [10]. Equating both linear functions results in a product of pure moment $M$ and radius of curvature $\rho$

$$M\rho = EI\,, \tag{6.1}$$

which is constant, since the bending rigidity expression $EI$ on the right side of (6.1) depends on both constant cross-section and material parameters. The bending energy $U$ is elastic strain energy stored in the beam of length $L$ [11].

$$U = \frac{1}{2}\int_0^L \frac{M}{\rho}\,ds$$

Using relation (6.1) and its constancy herein, then integrating simplifies to

$$U = \frac{EI}{2}\int_0^L \frac{1}{\rho^2}\,ds = \frac{EI}{2}\frac{L}{\rho^2} \tag{6.2}$$

and thus leads to

**Theorem 4**:

Let $EI$ be the arc's rigidity expression, $\theta$ be its enclosed angle and $c$ be the distance from its startpoint to its endpoint. The bending energy is then

$$U = EI \frac{\theta \sin \frac{\theta}{2}}{c} \tag{6}$$

which is always positive.

*Proof.* Substituting $\rho = R$ in equation (6.2) and with $L = \theta R$ according to ,(4) and $R$ from equation (1) results in relation (6). The positive sign follows from the same argumentation as in theorem 2. □

Subsequently we will normalize the beam's bending rigidity to $EI = 1$ as a permissible simplification.

## 3. Polyarcs and Arc Splines

While a plane polygonial curve can be described as a sequence of $(x_i, y_i)$ pairs, a polyarc is formulated as a sequence of $(x_i, y_i, \theta_i)$ triples, where arc angle $\theta_i$ belongs to circular arc from point $i$ to point $i+1$.

```
[ {x1, y1, θ1},
  {x2, y2, θ2},
  {x3, y3, θ3},
        ⋮
  {xn, yn, θn} ]
```

Table: Polyarc notation in JSON syntax.

Using this convention, any polygonial curve can be seen as a polyarc with all $\theta_i$ set to zero. The arcs of a polyarc do not depend on neighbour arc geometry (Fig. 4a). Both polygonial curves and polyarcs have $G^0$ continuity. So two adjacent arcs only have their join point in common

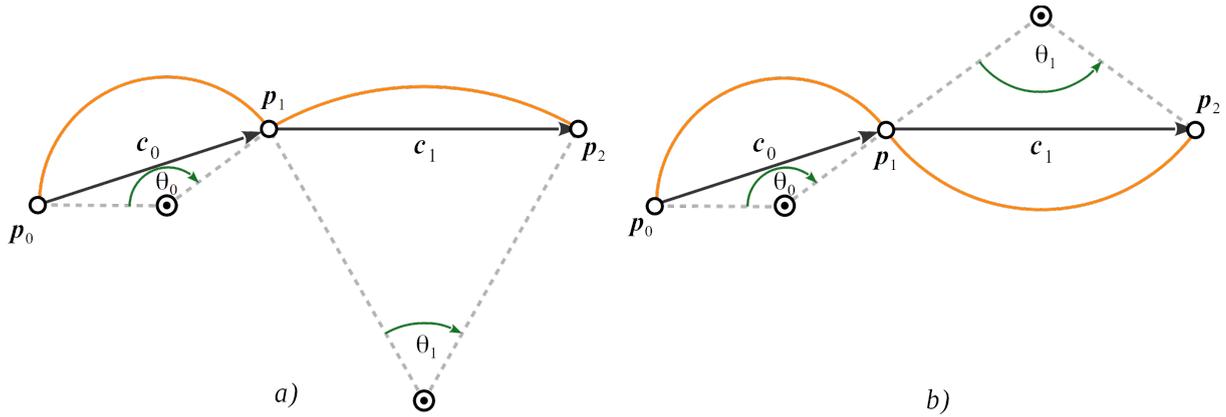

Fig. 4: Polyarc and arc-spline.

In contrast, an arc spline is a curve with $G^1$ continuity, where two adjacent arcs have a common tangent in their join point (Fig. 4b). A straight line through the centers of two adjacent arcs in an arc-spline always contains their join point.

Nevertheless, an arc spline is also always a polyarc – accidentally having the additional property of $G^1$ continuity.

## 3.1 Polyarc Properties

In contrast to cubic splines the total length of a polyarc can be not only approximated, but determined exactly by summing up the length of all individual arcs $L_i = s_i(1)$ according to equation (4).

$$L = \sum_i^n L_i$$

The area of a not necessarily closed polyarc can be calculated as the sum of the area of its polygonial curve and all individual arc segment areas according to equation (5), i.e.

$$A = \frac{1}{2} \sum_i^n \tilde{\mathbf{p}}_i \mathbf{p}_{i+1} + \sum_i^n A_{seg,i}$$

For closed polygons simply set $\mathbf{p}_{i+1} = \mathbf{p}_0$. Please also note, that we won't restrict its orientation (counter- or clockwise) and selfintersection. Simply interpret the result area with these aspects in mind.

We even can calculate the bending energy of the total curve via the sum of the bending energy of all individual arcs according to equation (6) as we want use it later.

$$U = \sum_i^n U_i \,,$$

## 3.2 Arc Spline Geometry

The start tangent and the end tangent of circular arc $i$ include equal but opposite angles with their base vector $\mathbf{c}_i$. Having given either angle $\theta_i$ or start vector $\mathbf{t}_i$ we are able to determine the respective other value.

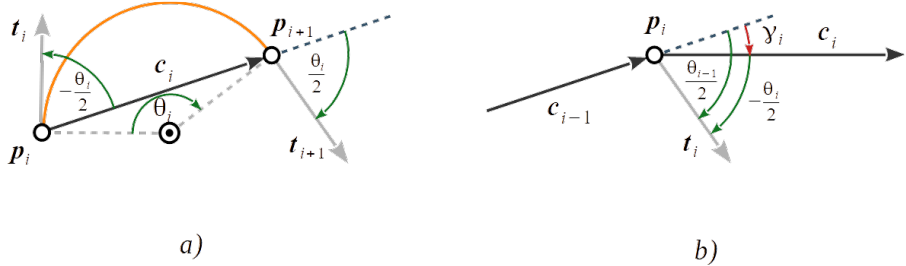

a)          b)

Fig. 5: Arc spline tangents and angles.

**Lemma 3**:
Let $\mathbf{c}_i$ be the edge vector from vertex $\mathbf{p}_i$ to vertex $\mathbf{p}_{i+1}$ of a polygonial curve. Let further be $\theta_i$ the signed angle of the circular arc erected over that edge vector. The unit vector $\mathbf{t}_i$ tangential to the arc in $\mathbf{p}_i$ is then

$$\mathbf{t}_i = \frac{\cos \frac{\theta_i}{2} \mathbf{c}_i - \sin \frac{\theta_i}{2} \tilde{\mathbf{c}}_i}{c_i} \tag{7}$$

*Proof.* We yield the tangent vector $\mathbf{t}_i$ in vertex $\mathbf{p}_i$ by rotating the unit edge vector $\frac{\mathbf{c}_i}{c_i}$ by $-\frac{\theta_i}{2}$ (Fig. 5a). □

**Lemma 4**:
Let $\mathbf{c}_i$ be the edge vector from vertex $\mathbf{p}_i$ to vertex $\mathbf{p}_{i+1}$ of a polygonial curve and $\mathbf{t}_i$ the unit vector tangential to the arc over $\mathbf{c}_i$ in $\mathbf{p}_i$. The arc angle $\theta_i$ is then given by

$$\tan \frac{\theta_i}{2} = \frac{\tilde{\mathbf{t}}_i \mathbf{c}_i}{\mathbf{t}_i \mathbf{c}_i} \tag{8}$$

*Proof.* Multiplying equation (7) by unit edge vector $\frac{\mathbf{c}_i}{c_i}$ and its orthogonal complement $\frac{\tilde{\mathbf{c}}_i}{c_i}$ respectively isolates the trigonometric terms in equation (7) and gives

$$\cos\frac{\theta_i}{2} = \frac{\mathbf{t}_i\,\mathbf{c}_i}{c_i} \quad and \quad \sin\frac{\theta_i}{2} = -\frac{\tilde{\mathbf{c}}_i\,\mathbf{t}_i}{c_i} = \frac{\tilde{\mathbf{t}}_i\,\mathbf{c}_i}{c_i}.$$

Equation (8) results by dividing the second term by the first. □

Focusing on the geometric relations in vertex $i$, where two succeeding edge vectors meet, gives us the dependency between adjacent arc angles in an arc spline. We recognize from Figure 5b the angular relationship

$$\frac{\theta_i}{2} + \frac{\theta_{i-1}}{2} - \gamma_i = 0 \quad where \quad \tan\gamma_i = \frac{\tilde{\mathbf{c}}_{i-1}\mathbf{c}_i}{\mathbf{c}_{i-1}\mathbf{c}_i} \tag{9.1}$$

Given that, we will be able to show, that any arc angle $\theta_i$ only depends on start angle $\theta_0$ and its preceding exterior angles $\gamma_j$ of the underlying polygonial curve.

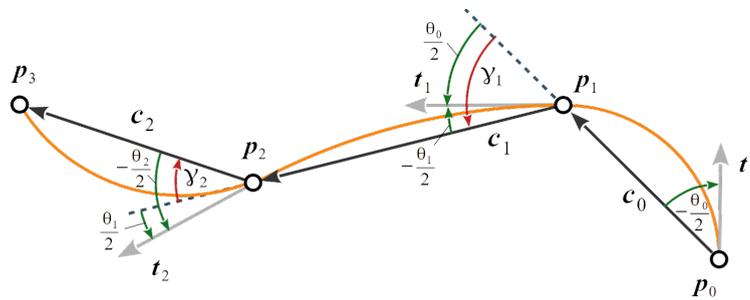

Fig. 6: Arc spline angles.

**Theorem 5**:
Let $\theta_0$ be the arc angle of the first arc of an arc spline, $\theta_i$ the angle of the $i^{th}$ following arc and let $\gamma_j$ be all preceding exterior angles. The value of arc angle $\theta_i$ then obeys

$$\theta_i = (-1)^i \left[\theta_0 + 2\sum_{j=1}^{i}(-1)^j\gamma_j\right]. \tag{9}$$

It depends on start angle $\theta_0$ and the alternating sum of preceding exterior angles of the polygonial curve.

*Proof.* Writing down the first three angular relations according to equation (9.1) reads (Fig. 6)

$$\frac{\theta_1}{2} + \frac{\theta_0}{2} - \gamma_1 = 0$$
$$\frac{\theta_2}{2} + \frac{\theta_1}{2} - \gamma_2 = 0$$
$$\frac{\theta_3}{2} + \frac{\theta_2}{2} - \gamma_3 = 0$$
$$\vdots$$

Resolving each equation for current $\theta$ and reintroducing the previous gives:

$$\theta_1 = -\theta_0 + 2\gamma_1 = -(\theta_0 + 2(-\gamma_1))$$
$$\theta_2 = -\theta_1 + 2\gamma_2 = \phantom{-}(\theta_0 + 2(-\gamma_1 + \gamma_2))$$
$$\theta_3 = -\theta_2 + 2\gamma_3 = -(\theta_0 + 2(-\gamma_1 + \gamma_2 - \gamma_3))$$
$$\vdots$$

Generalizing this sequence for angle $\theta_i$ leads to equation (9). □

> **Remark.** You should be aware of the fact, that an arc angle $\theta$ by equation (9) might get a value outside of the interval $(-2\pi, 2\pi)$. In that case one should take the complement angle $\theta \mp 4\pi$ instead.

## 3.2 One Parameter Family of Arc Splines

According to equation (9), the arc spline of a given polygonial curve – having $G^1$ continuity – is determined exclusively by the arc angle $\theta_0$ of the first segment. The result is a one parameter family of arc spline curves. The start arc angle $\theta_0$ in Figure 7 takes values from $-180°$ to $180°$ in $60°$ steps.

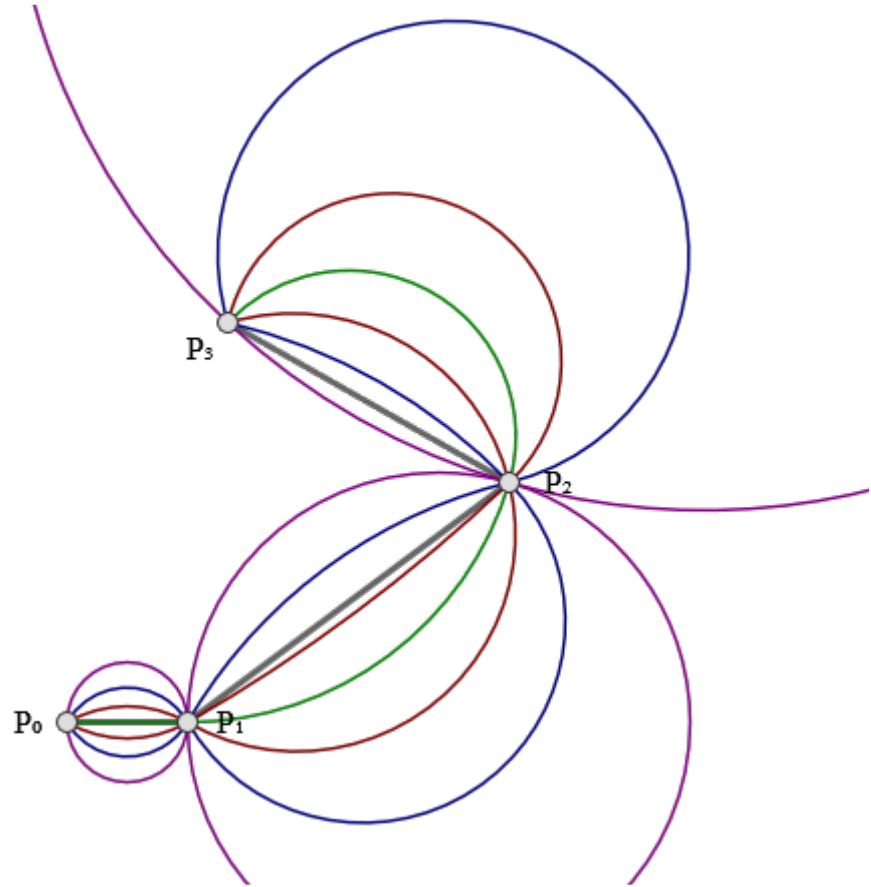

Fig. 7: Example of an one parameter family of arc splines.

> **Remark.** To be precise, it is not necessarily the first segment's angle, but the angle of any single segment, that must be specified. However, for the sake of simplicity, we will stick to that first angle.

Note that $G^1$ continuity is not generally guaranteed for the last segment of closed polygonial curves, but occurs only in special cases, such as symmetrical ones. If we want $G^1$ continuity also for closed curves, we might use a biarc for the closing segment, since in general two circular arcs are required to meet two adjacent points and their given tangents.

## 3.3 Smoothing of Arc Splines

Now having identified $\theta_0$ as a suitable single parameter of a family of arc splines belonging to a given polygonial curve, we want to discuss some criterions for a good choice of $\theta_0$. For this we take the approach of optimizing an arc spline property and discuss minimization of

- arc spline length.
- areal difference to polygonial curve.
- arc spline bending energy.

Whatever we choose from this list, we will perform a global minimization with it.

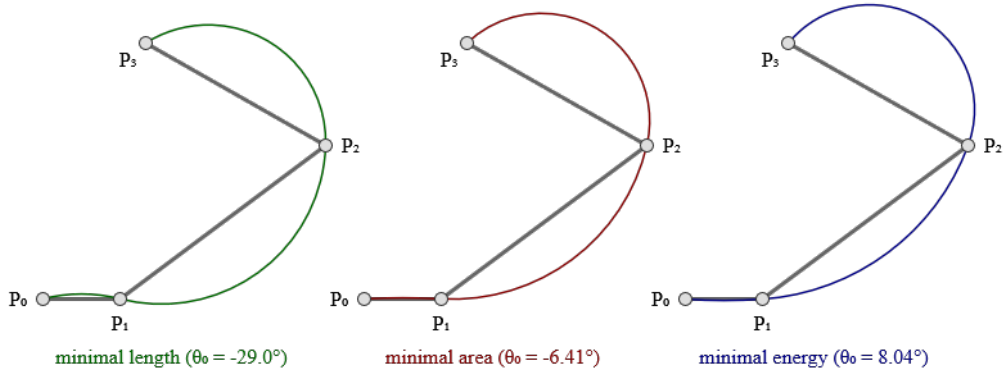

Fig. 8: Global minimization of different arc spline properties.

For minimizing the curve length, we need to find the minimum of the sum of all segment arc lengths according to equation (4), i.e. $L = \sum_i L_i$. In the case of minimizing the area we sum up the absolute value of all individual arc segment areas according to equation (5), $A = \sum_i |A_i|$. Minimizing the overall bending energy uses equation (6), i.e. $U = \sum_i U_i$, while the bending stiffness $EI$ is normalized to $1\,Nmm^2$.

In Table 1 are summerized the segment equations used, start angle found, length, area and bending energy for each of the three cases in Figure 8.

Table 1: Length, area and bending energy of example in Figure 8.

| minimize by | use equation | $\theta_0\,[°]$ | $L\,[mm]$ | $A\,[mm^2]$ | $U\,[Nmm]$ |
|---|---|---|---|---|---|
| Length | $L_i = \dfrac{c}{2}\dfrac{\theta}{\sin\frac{\theta}{2}}$ | $-29.0$ | 488 | 12420 | 596 |
| Area | $A_i = \dfrac{c^2}{4}\dfrac{\theta - \sin\theta}{1 - \cos\theta}$ | $-6.41$ | 493 | 12139 | 576 |
| Energy | $U_i = EI\,\dfrac{\theta \sin\frac{\theta}{2}}{c}$ | 8.04 | 501 | 12347 | 572 |

Calculating global extreme values of composite curves is usually not analytically feasible. Instead we apply the *golden-section search* as a numerical approach. The *golden-section search* is an efficient method to progressively reduce the interval locating the minimum or maximum. There are various implementations in different computer languages of this algorithm available [12].

```javascript
// recursive version of golden section search for minimum of
// function 'fn(x)' in interval '[lo,up]' using tolerance 'tol'.
function gss(fn,lo,up,tol) {
    const delta = up-lo;
    if (delta > tol) {
        const _lo = up - delta*invratio;
        const _up = lo + delta*invratio;
        return fn(_lo) < fn(_up) ? gss(fn,lo,_up,tol)
                                 : gss(fn,_lo,up,tol);
    }
    return (up+lo)/2;
}
const invratio = (Math.sqrt(5)-1)/2; // inverse of golden ratio
```

Listing 1: JavaScript function for recursive golden section search

Listing 1 shows a recursive implementation of the algorithm in JavaScript, which one might be able to port to another language quite easily. This algorithm proves to be highly robust. It finds each minimum of the examples in Figure 8 quite quickly in 15 steps within a generous starting range of $[-344°, 344°]$ and an angular tolerance of $0.6°$.

# 4. Conclusion

In this paper circular arcs are discussed from a symplectic geometry point of view.

From two different popular parameter sets – center based and endpoint based parameters – the latter is more appropriate for use with polyarcs and arc splines. In a polyarc the individual arcs are not related to each other ($G^0$ continuity) in contrast to an arc spline, where adjacent arcs share a common tangent in their join point ($G^1$ continuity). An arc spline is basically a special polyarc, which is build over an underlying polygonial curve by definition.

So with any polygonial curve we have a one parameter family of arc splines. Since we have chosen the first arc angle $\theta_0$ of the arc spline to be that parameter in question, we are able to define any subsequent arc angle $\theta_i$ by equation (9). The question for an optimal parameter leads to the search of global minima of curve length, bending energy or areal difference to the underlying polygonial curve.

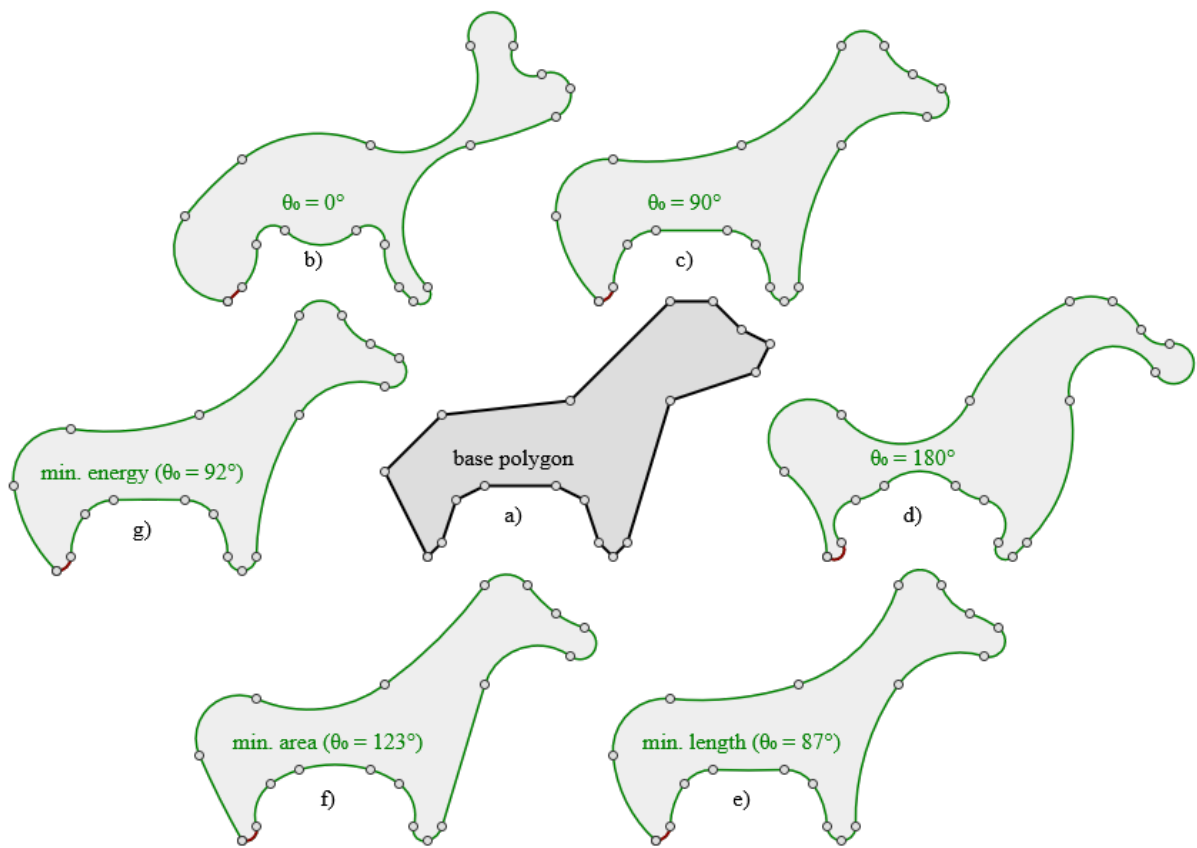

Fig. 9: Global minimization of different arc spline properties.

An illustrative example for an arc spline family generated from a base polygon (Fig.9a) is shown in Figure 9. The start arc angle $\theta_0$ – colored red at the bottom of the hind leg – has either explicitly given values (Fig.9b-d) or is numerically found by minimization criterions (Fig.9e-g).

As a result of this paper it has been shown one more time, that a consequent use of symplectic geometry in $\mathbb{R}^2$ may lead to new elegant result equations. Arc splines generated this way are a simple and compact alternative to biarcs, as long as no tangents at the vertices need to be taken into account.